\documentstyle[12pt]{article}
\begin{document}
\noindent
\begin{center}
{\large\bf FRACTIONAL REACTION-DIFFUSION EQUATIONS}\\[0.5cm]
R.K. SAXENA\\
Department of Mathematics and Statistics, Jai Narain Vyas University\\
Jodhpur, 342005, India\\[0.5cm]
A.M. MATHAI\\
Department of  Mathematics and Statistics, McGill University\\
Montreal, Canada  H3A 2K6\\[0.5cm]
H.J. HAUBOLD\\
Office for Outer Space Affairs, United Nations\\
P.O.Box 500, A–1400 Vienna, Austria\\[0.5cm]
\end{center}
\bigskip
\noindent
{\bf Abstract.}
In a series of papers, Saxena, Mathai, and Haubold (2002, 2004a, 2004b) derived solutions of a number of fractional kinetic
equations  in terms of generalized Mittag-Leffler functions which provide the extension of the work of Haubold and 
Mathai (1995, 2000). The subject of the present paper is to investigate the solution of a fractional reaction-diffusion equation. The results derived are of general nature and  include the results reported earlier by many authors, notably by  Jespersen, Metzler, and  Fogedby (1999) for anomalous diffusion  and del-Castillo-Negrete, Carreras, and Lynch (2003) for reaction-diffusion systems with L\'evy flights. The solution has been developed in terms of the H-function in a compact form with the help of Laplace and Fourier transforms. Most of the results obtained are in a form suitable for numerical computation.
\section {Introduction} 
Reaction-diffusion  models have found numerous applications in pattern formation in biology, chemistry, and physics, 
see Murray (2003), Kuramoto (2003), Wilhelmsson and Lazzaro (2001), and Hundsdorfer and Verwer (2003). These systems indicate that diffusion can produce spontaneous formation of 
spatio-temporal patterns. For details, one can refer to the work of Nicolis and Prigogine (1977) and Haken (2004). A general model for reaction-diffusion 
systems is investigated by Henry and Wearne (2000, 2002) and Henry, Langlands, and Wearne (2005). 

The simplest reaction-diffusion models are of the form
\begin{equation} 
\frac{\partial N}{\partial t}=d\frac{\partial^2N}{\partial x^2}+F(N), N=N(x,t),
\end{equation} 
where  $d$ is the diffusion coefficient and $F(N)$ is a nonlinear function representing  reaction kinetics. It is interesting 
to observe that for $F(N) =\gamma  N(1-N^2),$ eq. (1) reduces to the Fisher-Kolmogorov  equation and if we set  
$F(N) = N(1- N^2),$ it gives rise to the real Ginsburg-Landau equation. Del-Castillo-Negrete, Carreras, and Lynch (2002) 
studied the front propagation and segregation in a system of reaction-diffusion equations with cross-diffusion. 
Recently, del-Castillo-Negrete, Carreras, and Lynch (2003) discussed the dynamics in reaction-diffusion systems with non-Gaussian 
diffusion caused by asymmetric L\'evy flights and solved the following model
\begin{equation}
\frac{\partial N}{\partial t}=\eta D_x^\alpha N+F(N), N=N(x,t),
\end{equation}
with $F=0$.

In this paper we present a solution of a more general  model of reaction-diffusion systems (2) in which 
$\frac{\partial N}{\partial t}$
has been replaced by $\frac{\partial^\beta N}{\partial t^\beta},\beta > 0.$ This new model  extends the work of 
Jespersen, Metzler, and Fogedby (1999) and del-Castillo-Negrete, Carreras, and Lynch (2003). Most of the results are obtained in a compact form suitable for numerical computation.

A generalization of the Mittag-Leffler  function (Mittag-Leffler, 1903, 1905)
\begin{equation}
E_\alpha(z):=\sum^\infty_{n=0}\frac{z^n}{\Gamma(n\alpha+1)}, \alpha \in C, Re(\alpha)>0,
\end{equation}
was introduced by Wiman (1905)  in the generalized form 
\begin{equation}
E_{\alpha,\beta}(z):= \sum^\infty_{n=0}\frac{z^n}{\Gamma(n\alpha+\beta)},\alpha, \beta \in C, Re(\alpha)>0.
\end{equation}
The main results of these functions are available in the handbook of Erd\'elyi, Magnus, Oberhettinger, and Tricomi 
(1955, Section 18.1) and monographs by Dzherbashyan (1966, 1993).  
    
The H-function is defined by means of a Mellin-Barnes type integral in the following manner 
(Mathai and Saxena, 1978)
\begin{eqnarray}
H^{m,n}_{p,q}(z)& = &H^{m,n}_{p,q}\left[z\left|^{(a_p, A_p)}_{(b_q, B_q)}\right.\right]\\\nonumber 
&=& H^{m,n}_{p,q}\left[z\left|^{(a_1,A_1),\ldots,(a_p,A_p)}_{(b_1,B_1),\ldots,(B_q,B_q)}\right.\right]
=\frac{1}{2\pi i}\int_\Omega \Theta(\xi)z^{-\xi}d\xi, 
\end{eqnarray}              
where $i=(-1)^{1/2},$
\begin{equation}
\Theta(\xi)=\frac{\left[\Pi^m_{j=1}\Gamma(b_j+B_j\xi)\right]\left[\Pi^n_{j=1}\Gamma(1-a_j-A_j\xi)\right]}
{\left[\Pi^q_{j=m+1}\Gamma(1-b_j-B_j\xi)\right]\left[\Pi^p_{j=n+1}\Gamma(a_j+A_j\xi)\right]},
\end{equation}
and an empty product is always interpreted as unity; $m, n, p, q \in N_0$ with  
$0\leq n \leq p, 1\leq m\leq q, A_i,B_j\in R_+,a_i,b_j \in R$ or $C(i=1,\ldots,p; j=1,\ldots,q)$ such that
\begin{equation}
A_i(b_j+k)\neq B_j(a_i-l-1), k,l \in N_0; i=1,\ldots,n; j=1,\ldots,m,
\end{equation}
where we employ the usual notations:$N_0=(0,1,2\ldots), R = (-\infty, \infty), R_+ =(0,\infty)$ and $C$ being the complex number field.  The contour $\Omega$  is either $L_{-\infty}, L_{+\infty},$ or $L_{i\gamma\infty}.$ These contours are defined  in the monographs  by Prudnikov, Brychkov, and Marichev (1989), Mathai (1993), and Kilbas and Saigo (2004). However, the  explicit definitions of these contours are given below to complete the definition of the H-function.\\

(i)  $\Omega = L_{-\infty}$ is a left loop situated in a horizontal strip starting at the point $-\infty+i\varphi-1$
and terminating at the point $-\infty+i\varphi_2$ with $-\infty<\varphi_1<\varphi_2<+\infty;$\\
 
(ii) $\Omega=L_{+\infty}$ is a right loop situated in a horizontal strip starting at the point $+\infty+i\varphi_1$ and terminating at the point $+\infty+i\varphi_2$  with $-\infty<\varphi_1<\varphi_2<+\infty.$\\

(iii) $\Omega=L_{i\gamma\infty}$ is a contour starting at the point $\gamma-i\infty$  and terminating at the point  $\gamma+i\infty$, where $\gamma\in R=(-\infty,+\infty)$.\\

A detailed and comprehensive account of the H-function is available from the monograph by Mathai and Saxena (1978), Prudnikov, Brychkov, and Marichev (1989), and Kilbas and Saigo (2004). The relation connecting $_p\Psi_q(z)$ and the H-function is given for the first time in the monograph by Mathai and Saxena (1978, p.11, Eq.1.7.8) as
\begin{equation}
_p\Psi_q\left[^{(a_1,A_1),\ldots,(a_p,A_p)}_{(b_1,B_1),\ldots,(b_q,B_q)}\left|z\right.\right]=H^{1,p}_{p,q+1}\left[-z\left|^{(1-a_1,A_1),\ldots,(1-a_p,A_p)}_{(0,1),(1-b_1,B_1),\ldots,(1-b_q,B_q)}\right.\right],
\end{equation}
where $_p\Psi_q(z)$ is Wright's generalized hypergeometric function (Wright,\\
1940); also see (Erd\'elyi, Magnus, Oberhettinger, and Tricomi, 1953, Section 4.1), defined by means of the series representation in the form
\begin{equation}                                         
_p\Psi_q(z)=\;_p\Psi_q\left[^{(a_p,A_p)}_{(b_q,B_q)}\left|z\right.\right]=\sum^\infty_{r=0}\frac{[\Pi^p_{j=1}\Gamma(a_j+A_jr)]z^r}{[\Pi^q_{j=1}\Gamma(b_j+B_jr)(r)!]},
\end{equation}
where $z\in C, a_i,b_j\in C,A_i, B_j \in R=(-\infty, \infty), A_i, B_j\neq 0\; (i=1,\ldots,p; j=1, \ldots,q),\sum^q_{j=1}B_j-\sum^p_{j=1}A_j>-1; C$ being the set of complex numbers and $\Gamma(z)$ is Euler's gamma function. 
This function  includes many special functions. It is interesting to observe that for $A_i=B_j=1, \forall i$ and $j$, eq. (9) reduces to a generalized hypergeometric function $_pF_q(z)$ as
\begin{equation}
_p\Psi_q\left.\left[^{(a_p,1)}_{(b_q,1)}\right|z\right]=\frac{\Pi^p_{j=1}\Gamma(a_j)}{\Pi^q_{j=1}\Gamma(b_j)}\;_pF_q(a_1,\ldots,a_p;b_1,\ldots,b_q;z),
\end{equation}
where $a_j\neq-\nu(j=1,\ldots,p$ and $\nu=0,1,2,\ldots.);p<q$ or $p=q,|z|< 1.$
Prior to (9), Wright (1933) introduced a special case of (9) in the form 
$$\Phi(a,b;z)=\;_0\Psi_1\left.\left[^{-}_{(b,a)}\right|z\right]=\sum^\infty_{r=0}\frac{1}{\Gamma(ar+b)}\frac{z^r}{(r)!},$$ 
which widely occurs in problems of fractional  diffusion. It has been shown by Saxena, Mathai, and Haubold (2004b) that 
\begin{eqnarray}
E_{\alpha,\beta}(z)&=&_1\Psi_1\left.\left[^{(1,1)}_{(\beta,\alpha)}\right|z\right]\\                    
&=& H^{1,1}_{1,2}\left[-z\left|^{(0,1)}_{(0,1),(1-\beta,\alpha)}\right.\right].
\end{eqnarray}
If we further take $\beta=1$  in (11) and (12), we find that
\begin{eqnarray}
E_{\alpha,1}(z) &=&E_\alpha(z)=\;_1\Psi_1\left.\left[^{(1,1)}_{(1,\alpha)}\right|z\right]\\
&=& H^{1,1}_{1,2}\left[-z\left|^{(0,1)}_{(0,1),(0,\alpha)}\right.\right],
\end{eqnarray}                             
where $Re(\alpha)>0, \alpha\in C.$ 

From Mathai and Saxena (1978)  and Prudnikov, Brychkov, and Marichev
(1989, p.355, Eq.2.25.3), it follows that the Laplace transform of the H-function is given by
\begin{equation}
L\left\{t^{\rho-1}H^{m,n}_{p,q}\left[zt^\sigma\left|^{(a_p,A_p)}_{(b_q,B_q)}\right.\right]\right\}= s^{-\rho}H^{m,n+1}_{p+1,q}\left[zs^{-\sigma}\left|^{(1-\rho,\sigma),(a_p,A_p)}_{(b_q,B_q)}\right.\right],
\end{equation}
where $\sigma>0, Re(s)>0, Re[\rho+\sigma^{min}_{1\leq j\leq m}(\frac{b_j}{B_j})]>0, |arg z|<[\pi/2]\theta, \theta>0;$
$\theta=\sum^n_{j=1}A_j-\sum^p_{j=n+1}A_j+\sum^m_{j=1}B_j-\sum^q_{j=m+1}B_j.$

By virtue of the cancellation law for the H-function (Mathai and Saxena, 1978), it can be readily seen that 
\begin{equation}
L^{-1}\left\{s^{-\rho}H^{m,n}_{p,q}\left[zs^\sigma\left|^{(a_p,A_p)}_{(b_q,B_q)}\right.\right]\right\}=t^{\rho-1}H^{m,n}_{p+1,q}\left[zt^{-\sigma}\left|^{(a_p,A_p),(\rho,\sigma)}_{(b_q,B_q)}\right.\right],
\end{equation}
where $\sigma >0 ,Re(s)> 0, Re[\rho+\sigma^{max}_{1\leq j\leq n}(\frac{1-a_j}{A_j})]>0, |arg z| < \frac{1}{2}\pi\theta_1,\theta_1>0;$\\
$\theta=\theta-a$. Two interesting special cases of (16) are worth mentioning. If we employ the identity (Mathai and Saxena, 1978)
\begin{equation}
H^{1,0}_{0,1}\left[x\left|_{(\alpha,1)}\right.\right]=x^\alpha exp(-x),
\end{equation}
we obtain 
\begin{equation}
L^{-1}[s^{-\rho}exp(-zs^\sigma)]=t^{\rho-1}H^{1,0}_{1,1}\left[zt^{-\sigma}\left|^{(\rho,\sigma)}_{(0,1)}\right.\right],
\end{equation}
where $ Re(s)>0, \sigma>0.$\\
Further if we use the identity (Mathai and Saxena, 1978)
\begin{equation}
H^{2,0}_{0,2}\left[x\left|_{(\frac{\nu}{2},1)(-\frac{\nu}{2},1)}\right.\right]=2K_\nu(2x^{1/2}),
\end{equation}
eq.(16) yields
\begin{equation}
2L^{-1}[s^{-\rho}K_\nu(zs^\sigma)]=t^{\rho-1}H^{2,0}_{1,2}\left[\frac{z^2t^{-2\sigma}}{4}\left|^{(\rho,2\sigma)}_{(\frac{\nu}{2},1),(-\frac{\nu}{2},1)}\right.\right],
\end{equation}     
where  $Re(\rho)> 0, Re(z^2)>0, Re(s)>0,$ and $K_\nu(.)$  is the Bessel function of the third kind .
                      
    In view of the result of Saxena,  Mathai, and Haubold (2004a, p.49), also see
Prudnikov, Brychkov, and Marichev (1989, p.355, eq.(2.25.3.2)), the cosine transform of the H-function is given by 
\begin{eqnarray} 
&&\int^\infty_0 t^{\rho-1}cos(kt)H^{m,n}_{p,q}\left[at^\mu\left|^{(a_p,A_p)}_{(b_q,B_q)}\right.\right]dt\\\nonumber  
&=&\frac{\pi}{k^\rho}H^{n+1,m}_{q+1,p+2}\left[\frac{k^\mu}{a}\left|^{(1-b_q,B_q),(\frac{1+\rho}{2},\frac{\mu}{2})}_{(\rho,\mu),(1-a_p,A_p),(\frac{1+\rho}{2},\frac{\mu}{2})}\right.\right],
\end{eqnarray}
where $Re[\rho+\mu^{min}_{1\leq j\leq m}]>1, |arg\; a|<\frac{1}{2}\pi\theta; \theta>0,\theta$ is defined with the result eq. (15).

The Riemann-Liouville fractional integral of order $\nu$ is defined by 
(Miller and Ross, 1993, p.45)
\begin{equation}
_0D_t^{-\nu}f(t)=\frac{1}{\Gamma(\nu)}\int_0^t(t-u)^{\nu-1}f(u)du,
\end{equation} 
where  $Re(\nu) > 0.$

    Following  Samko, Kilbas, and Marichev (1990, p.37), we define the fractional derivative  for $\alpha>0$ in the form 
\begin{equation}
_0D_t^\alpha f(t)=\frac{1}{\Gamma(n-\alpha)}\frac{d^n}{dt^n}\int_0^t\frac{f(u)du}{(t-u)^{\alpha-n+1}},\;\;n=[\alpha]+1,
\end{equation}      
where  $[\alpha]$   means the integral part of the number $\alpha$.\\
In particular, if  $0<\alpha<1,$
\begin{equation}
_0D_t^\alpha f(t)=\frac{d}{dt}\frac{1}{\Gamma(1-\alpha)}\int_0^t\frac{f(u)du}{(t-u)^\alpha},
\end{equation}
and if $\alpha=n\;\;\in N= \left\{1,2,\ldots,\right\},$ then
\begin{equation}
_0D_t^nf(t)=D^n f(t) (D=d/dt),
\end{equation}
is the usual derivative of order $n$.

     From  Erd\'elyi, Magnus, Oberhettinger, and Tricomi (1954b, p. 182), we have 
\begin{equation}
L\left\{_0D_t^{-\nu}f(t)\right\}=s^{-\nu}F(s),
\end{equation}
where
\begin{equation}
F(s)=L\left\{f(t);s\right\}=f^*(s)=\int^\infty_0 exp(-st)f(t)dt, Re(s)>0.
\end{equation}
The Laplace transform of the fractional derivative  is given by  
Oldham and Spanier (1974, p.134, eq.(8.1.3))
\begin{equation}
L\left\{_0D_t^\alpha f(t)\right\}=s^\alpha F(s)-\sum^n_{r=1}s^{r-1}\;_0D_t^{\alpha-r}f(t)|_{t=0}.
\end{equation}
In certain boundary-value problems, the following fractional derivative of order $\alpha>0$ is introduced by Caputo (1969) in the form
\begin{eqnarray}
D_t^\alpha f(t) & = & \frac{1}{\Gamma(m-\alpha)}\int_0^t\frac{f^{(m)}(\tau)d\tau}{(t-\tau)^{\alpha+1-m}},\\
&&m-1<\alpha\leq m, Re(\alpha)>0, m\in N.\nonumber\\
& = & \frac{d^mf}{dt^m}, \mbox{if}\;\; \alpha=m.
\end{eqnarray}
                                                                       
Caputo (1969)  has  given the Laplace transform of this derivative as 
\begin{equation}
L\left\{D_t^\alpha f(t); s\right\}= s^\alpha F(s)-\sum^{m-1}_{r=0}s^{\alpha-r-1}f^r(0+), m-1<a\leq m.
\end{equation}
The above formula is very useful in deriving the solution of differ-integral equations of fractional order governing certain physical problems of reaction and diffusion.
We also need the Weyl fractional operator defined by
\begin{equation}
_{-\infty}D_x^\mu f(t)=\frac{1}{\Gamma(n-\mu)}\frac{d^n}{dt^n}\int_{-\infty}^t\frac{f(u)du}{(t-u)^{\mu-n+1}},
\end{equation}                                                 
where  $n = [\mu]$  is an integral  part of $\mu > 0.$\\
Its Fourier transform is (Metzler and Klafter, 2000, p.59, A.11, 2004)
\begin{equation}
F\left\{_{-\infty}D_x^\mu f(x)\right\}=(ik)^\mu \tilde{f}(k),
\end{equation}
where we define the Fourier transform as 
\begin{equation}
\tilde{h}(q)=\int^\infty_{-\infty}h(x)\mbox{exp}(iqx)dx.
\end{equation}
Following the convention initiated by Compte (1996), we suppress the imaginary unit in Fourier space by adopting the slightly modified form of above result in our investigations (Metzler and Klafter, 2000, p.59, A.12, 2004) 
\begin{equation}
F\left\{_{-\infty}D^\mu_x f(x)\right\}=-|k|^\mu \tilde{f}(k)       
\end{equation}
instead of (33). Finally we also need the following property of the H-function 
(Mathai and Saxena, 1978)
\begin{equation}
H^{m,n}_{p.q}\left[x^\delta\left|^{(a_p,A_p)}_{(b_q,B_q)}\right.\right]= \frac{1}{\delta}H^{m,n}_{p,q}\left[x\left|^{(a_p,A_p/\delta)}_{(b_q,B_q/\delta)}\right.\right],
\end{equation}
where $\delta >0.$

\section{The Fractional Reaction-Diffusion Equation}
In this section we will investigate the solution of the reaction-diffusion equation (37). The result is given in the form of the following
theorem.\\ 
{\bf Theorem.} Consider the following fractional reaction-diffusion model 
\begin{equation}
\frac{\partial^\beta N(x,t)}{\partial t^\beta}=\eta\;\; _{-\infty}D^\alpha_x N(x,t)+\varphi(x,t); \eta,t>0, x\in R, 0<\beta \leq 2,
\end{equation}
with the initial condition 
\begin{equation}
N(x,0) = f(x), N_t(x,0)=g(x)\; \mbox{for}\;\; x\in R,
\end{equation}          
where $N_t(x,0)$ means the first partial derivative of $N(x,t)$ with respect to $\varphi$ evaluated at $t=0, \eta$ is a diffusion constant and $\varphi(x,t)$  is a nonlinear function belonging to the area of reaction-diffusion. Then for the solution of (37), subject to the initial conditions (38), there holds the formula

\begin{eqnarray}
N(x,t)&=&\frac{1}{2\pi}\int^\infty_{-\infty} \tilde{f}(k)E_{\beta,1}(-\eta|k|^\alpha t^\beta)exp(-ikx)dk\\\nonumber
&+& \frac{1}{2\pi}\int^\infty_{-\infty}t \tilde{g}(k)E_{\beta,2}(\eta|k|^\alpha t^\beta)exp(-ikx)dk\\\nonumber
&+&\frac{1}{2\pi}\int^t_0\xi^{\beta-1}\int^\infty_{-\infty}\tilde{\varphi}(k,t-\xi)E_{\beta,\beta}(-\eta|k|^\alpha \xi^\beta)\mbox{exp}(-ikx)dk d\xi.
\end{eqnarray}
{\bf Proof.} If we apply the Laplace transform with respect to the time variable t and Fourier transform with respect to space variable $x$ and use the initial conditions and (32), then the given equation transforms into the form
\begin{equation}
N^{\tilde{*}}(k,s)=\frac{\tilde{f}(k)s^{\beta-1}}{s^\beta+\eta|k|^\alpha}+\frac{\tilde{g}(k)s^{\beta-2}}{s^\beta +\eta|k|^\alpha}+\frac{\varphi^{\tilde{*}}(k)}{s^\beta+\eta|k|^\alpha}.
\end{equation}
On taking the inverse Laplace transform of (40) and applying the result 
\begin{equation}
L^{-1}\left\{\frac{s^{a-1}}{a+s^\beta}\right\}=t^{\beta-\alpha}E_{\beta,\beta-\alpha+1}(-at^\beta),
\end{equation}
where $Re(s)> 0, Re(\beta-\alpha+1)>0$, it is seen that
\begin{eqnarray}
\tilde{N}(k,t)&=& \tilde{f}(k)E_{\beta,1}(-\eta|k|^\alpha t^\beta)+\tilde{g}(k)tE_{\beta,2}(-\eta|k|^\alpha t^\beta)\\\nonumber
&+&\int_0^t\tilde{\varphi}(k,t-\xi)\xi^{\beta-1}E_{\beta,\beta}(-\eta|k|^\alpha \xi^\beta)d\xi.
\end{eqnarray} 
The required solution (39) is now obtained by taking the inverse Fourier transform of (38). Thus, we have
\begin{eqnarray*}
N(x,t)&=&\frac{1}{2\pi}\int^\infty_{-\infty}\tilde{f}(k)E_{\beta,1}(-\eta|k|^\alpha t^\beta)\mbox{exp}(-ikx)dk\\
&+& \frac{1}{2\pi}\int^\infty_{_\infty}t\tilde{g}(k)E_{\beta,2}(-\eta|k|^\alpha t^\beta)\mbox{exp}(-ikx)dk\\
&+&\frac{1}{2\pi}\int^t_0\xi^{\beta-1}\int_{-\infty}^\infty\varphi^*(k,t-\xi)E_{\beta,\beta}(-\eta|k|^\alpha \xi^\beta)\mbox{exp}(-ikx)dk d\xi.
\end{eqnarray*}
This completes the proof of the theorem.
If we set $\alpha =2$, we obtain the result given by Debnath (2003).\\
{\bf Note.} By virtue of the identity (12), the solution (39) can be expressed in terms of the H-function as can be seen from the solutions given in the special cases of the theorem in the next section.

\section{Special Cases}

When $g(x) = 0,$ then applying the convolution theorem of the Fourier transform to the solution (35), the theorem yields\\
{\bf Corollary 1.1.} The solution of fractional reaction-diffusion equation 
\begin{equation}
\frac{\partial^\beta}{\partial t^\beta} N(x,t)-\eta _{-\infty} D^\alpha_x N(x,t)=\varphi(x,t), \;x \in R, t>0,\eta>0,
\end{equation}       
subject to the initial conditions
\begin{equation}
N(x,0)=f(x), N_t(x,0)=0\;\; \mbox{for}\; x \in R, 1<\beta\leq 2,
\end{equation}
where $\eta$ is  a diffusion constant and $\varphi(x,t)$ is a nonlinear function belonging to the area of reaction-diffusion, is given by 
\begin{eqnarray}
N(x,t) & = & \int^\infty_{-\infty} G_1(x-\tau, t) f(\tau)d\tau\\\nonumber
& + & \int^t_0(t-\xi)^{\beta-1}\;\;\int^x_0G_2(x-\tau, t-\xi)\varphi(\tau,\xi)d\tau d\xi,
\end{eqnarray}
where
\begin{eqnarray}
G_1(x,t) & = &\frac{1}{2\pi}\int^\infty_{-\infty}exp(-ikx)E_{\beta,1}(-\eta|k|^\alpha t^\beta)dk\\\nonumber
& = & \frac{1}{\pi \alpha}\int^\infty_0 cos(kx)H^{1,1}_{1,2}\left[k\eta^{1/\alpha}t^{\beta/\alpha}\left|^{(0,1/\alpha)}_{(0,1/\alpha),(0,\beta/\alpha)}\right.\right]dk\\\nonumber
& = & \frac{1}{\alpha|x|}H^{2,1}_{3,3}\left[\frac{|x|}{\eta^{1/\alpha}t^{\beta/\alpha}}\left|^{(1,1/\alpha),(1,\beta/\alpha),(1,1/2)}_{(1,1),(1,1/\alpha),(1,1/2)}\right.\right], \alpha>0,
\end{eqnarray}
\begin{eqnarray}
G_2(x,t) & = &\frac{1}{2\pi}\int^\infty_{-\infty}\mbox{exp}(-ikx)E_{\beta,\beta}(-\eta|k|^\alpha t^\beta)dk\\\nonumber
& = &\frac{1}{\pi\alpha}\int^\infty_0 cos(kx)H^{1,1}_{1,2}\left[k\eta^{1/\alpha}t^{\beta/\alpha}\left|^{(0,1/\alpha)}_{(0,1/\alpha),(1-\beta, \beta/\alpha)}\right.\right]dk\\\nonumber
& = &\frac{1}{\alpha|x|}H^{2,1}_{3,3}\left[\frac{|x|}{\eta^{1/\alpha}t^{\beta/\alpha}}\left|^{(1,1/\alpha),(\beta,\beta/\alpha),(1,1/2)}_{(1,1),(1,1/\alpha),(1,1/2)}\right.\right], \alpha>0.
\end{eqnarray}
If we set  $f(x)= \delta(x), \varphi\equiv 0, g(x)=0,$ where $\delta(x)$ is the Dirac-delta function, then we arrive at the following\\  
{\bf Corollary 1.2.} Consider the following reaction-diffusion model 
\begin{equation}
\frac{\partial^\beta N(x,t)}{\partial t^\beta}=\eta\;\;_{-\infty}D^\alpha_x N(x,t), \eta>0, x \in R, 0<\beta\leq 1,
\end{equation}       
with the initial condition  $N(x,t = 0)  = \delta(x),$ where $\eta$ is a diffusion constant  and $\delta(x)$ is the Dirac-delta function. Then  the solution of (44) is given by 
\begin{equation}
N(x,t)=\frac{1}{\alpha|x|}H^{2,1}_{3,3}\left[\frac{|x|}{(\eta t^\beta)^{1/\alpha}}\left|^{(1,1/\alpha),(1,\beta/\alpha),(1,1/2)}_{(1,1),(1,1/\alpha),(1,1/2)}\right.\right].
\end{equation}
In the case $\beta=1$, then in view  of the cancellation law for the H-function (Mathai and Saxena, 1978), (49) gives rise to the following result given by Jespersen, Metzler, and Fogedby (1999) and  recently by Del-Castillo-Negrete, Carreras, and Lynch (2003) in an entirely different form.\\

For the solution of fractional reaction-diffusion equation 
\begin{equation}
\frac{\partial}{\partial t}N(x,t)=\eta\;_{-\infty}D^\alpha_x N(x,t),
\end{equation}
with initial condition 
\begin{equation}
N(x,t=0)=\delta(x),
\end{equation}
there holds  the relation  
\begin{equation}
N(x,t)=\frac{1}{\alpha|x|}H^{1,1}_{2,2}\left[\frac{|x|}{\eta^{1/\alpha}t^{1/\alpha}}\left|^{(1,1/\alpha),(1,1/2)}_{(1,1),(1,1/2)}\right.\right],
\end{equation}
where $\alpha>0$.

It may be noted that (52) is a closed-form representation of a L\'evy stable law, see Metzler and Klafter (2000, 2004). It is interesting to note that as $\alpha\rightarrow 2$, the classical Gaussian solution is recovered  as 
\begin{eqnarray}
N(x,t) & = &\frac{1}{2|x|} H^{1,1}_{2,2}\left[\frac{|x|}{(\eta t)^{1/2}}\left|^{(1,1/2),(1,1/2)}_{(1,1),(1,1/2)}\right.\right]\nonumber\\
& = & \frac{1}{2|x|}H^{1,0}_{1,1}\left[\frac{|x|}{(\eta t)^{1/2}}\left|^{(1,1/2)}_{(1,1)}\right.\right]\\
& = & (4 \pi \eta t)^{-1/2} \mbox{exp}[-\frac{|x|^2}{4\eta t}].
\end{eqnarray}

It is useful to study the solution (49) due to its occurrence in certain fractional and diffusion models.  Now we proceed to find the fractional order moments of (49).\\
Here we remark that applying Fourier transform with respect to $x$ in (48) it is found that 
$$\frac{\partial^\beta}{\partial t^\beta}\tilde{N}(k,t)=-\eta|k|^\alpha \tilde{N}(k,t),$$
which is the generalized Fourier transformed diffusion equation, since for $\alpha=2$ and  for $\beta=1$, it reduces to Fourier transformed diffusion equation 
$$\frac{\partial \tilde{N}(k,t)}{\partial t}=\eta|k|^2 \tilde{N}(k,t),$$
being a diffusion equation, for a fixed wave number $k$ (Metzler and Klafter, 2000, 2004). Here $\tilde{N}(x,t)$  is the Fourier transform with respect to $x$ of $N(x,t).$

\section{Fractional Order Moments}
We now calculate the fractional order moments defined by
\begin{equation}
<|x(t)|^\delta>=\int^\infty_0|x|^\delta N(x,t)dx.
\end{equation}
Using the definition of the Mellin transform   
\begin{equation}
M\left\{q(t);s\right\}=\int^\infty_0t^{s-1}q(t)dt,
\end{equation}                             
we find from (49) that 
\begin{equation}
<|x(t)|^\delta>=\int^\infty_{-\infty}|x|^\delta N(x,t)dx
\end{equation}
\begin{equation}
<|x(t)|^\delta>=\frac{2}{\alpha}\int^\infty_0 x^{\delta-1}H^{2,1}_{3,3}\left[\frac{|x|}{\eta^{1/\alpha}t^{\beta/\alpha}}|^{(1,1/\alpha),(1,\beta/\alpha),(1,1/2)}_{(1,1),(1,1/\alpha),(1,1/2)}\right]dx.
\end{equation}
   Applying the integral formula
\begin{equation}
\int^\infty_0 x^{\xi-1} H^{m,n}_{p,q}\left[ax\left|^{(a_p,A_p)}_{(b_q,B_q)}\right.\right]dx=a^{-\xi}\Theta(-\xi),
\end{equation}
where  $^{- min}_{1\leq j\leq m}Re(\frac{b_j}{B_j})<Re(\xi)< max_{1\leq j\leq n}Re(\frac{1-a_j}{A_j}), |arg\; a|<\frac{1}{2}\pi \theta, \theta>0,$
$\theta$ is defined  in (6) and $\Theta(-\xi)$ in the definition of the H-function (6).
We see that 
\begin{equation}
<|x(t)|^\delta>=\frac{2}{\alpha}(\eta t^\beta)\frac{\delta}{\alpha}\frac{\Gamma(-\frac{\delta}{\alpha})\Gamma(1+\delta)\Gamma(1+\frac{\delta}{\alpha})}{\Gamma(-\frac{\delta}{2})\Gamma(1+\frac{\beta\delta}{\alpha})\Gamma(1+\frac{\delta}{2})}.
\end{equation}
where $Re(\delta)>-1$ and $Re(\delta+\alpha) > 0.$\\

 Two interesting special cases of (60) are given below.\\
(i) \begin{equation}\;\;\mbox{As}\;\; \delta\rightarrow 0,\mbox{then using the result}\; \frac{1}{\Gamma(z)}\sim z\;\;  \mbox{for}\;\; z<<1,
\end{equation}
we find that 
\begin{equation}
^{\mbox{lim}}_{\delta\rightarrow 0} <|x(t)|^\delta>=1.
\end{equation}
(ii) $\mbox{When}\;\; \alpha= 2,$ the linear time dependence is 
\begin{equation}
^{\mbox{lim}}_{\delta\rightarrow 2, \alpha\rightarrow 2}<|x(t)|^\delta =\frac{2\eta t^\beta}{\Gamma(1+\beta)}.
\end{equation}                                                             

\section{Behavior of the Solution in Equation (49)}
Eq. (49) can be written in terms of a Mellin-Barnes type integral as
\begin{equation}
N(x,t)=\frac{1}{\alpha|x|}\frac{1}{2\pi\omega}\int_L\frac{\Gamma(\frac{s}{\alpha})\Gamma(1-s)\Gamma(1-\frac{s}{\alpha})}{\Gamma(1-\frac{s\beta}{\alpha})\Gamma(1-\frac{s}{2})\Gamma(s/2)}\left[\frac{|x|}{\eta^{1\alpha}t^{\beta/\alpha}}\right]^s  ds.
\end{equation}
Let us assume that the poles of the integrand are simple. Now evaluating the sum of residues in ascending powers of $|x|$ by calculating the residues at the poles of  $\Gamma(1-s)$ at the points $s =1+\nu\; (\nu=0,1,2,\ldots)$ and  $\Gamma(1-\frac{s}{\alpha})$  at the points 
$s =\alpha(1+\nu)\;(\nu=0,1,2,\ldots)$, it is found that the series expansion of the general solution (45) is given by
\begin{eqnarray}
N(x,t)&=&\frac{1}{\alpha\eta^{1/\alpha}t^{\beta/\alpha}}\sum^\infty_{\nu=0}\frac{(-1)^\nu}{(\nu)!}\frac{\Gamma[\frac{1+\nu}{\alpha}]\Gamma(1-\frac{1+\nu}{\alpha})}{\Gamma[1-\frac{\beta(1+\nu)}{\alpha}]\Gamma[\frac{1-\nu}{2}]\Gamma[\frac{1+\nu}{2}]}\left[\frac{|x|}{\eta^{1/\alpha}t^{\beta/\alpha}}\right]^\nu\nonumber\\
&+&\frac{|x|^{\alpha-1}}{\eta t^\beta}\sum^\infty_{\nu=0}\frac{(-1)^\nu\Gamma[1-\alpha(1+\nu)]}{\Gamma[1-\beta(1+\nu)]\Gamma[1-\frac{\alpha}{2}(1+\nu)]\Gamma[\frac{\alpha(1+\nu)}{2}]}\left[\frac{|x|^\alpha}{\eta t^\beta}\right]^\nu
\end{eqnarray}

where $0<Re(\nu) < \alpha,\left\{\frac{|x|}{\eta t^\beta}\right\}<1.$
From (65), we infer that 
\begin{equation}
N(x,t)\sim  A + B |x|^{\alpha-1},\;\;\mbox{as}\;\; x\rightarrow 0, 
\end{equation}
where $A$ and $B$ are numerical constants.    
Further from the series expansion (61), it can be seen that the initial behavior is given by
\begin{eqnarray}
N(x,t)& \sim & \frac{\Gamma(\frac{1}{\alpha})\Gamma(1-\frac{1}{\alpha})}{\pi\alpha\eta^{1/\alpha}t^{\beta/\alpha}\Gamma(1-\frac{\beta}{\alpha})},\;\;\mbox{for}\; 1<\alpha< 2,\\
& \sim &\frac{\Gamma(\frac{1-\alpha}{2})|x|^{\alpha-1}}{\pi^{1/2}\eta 2^\alpha t^\beta\Gamma(1-\beta)\Gamma(\frac{\alpha}{2})},\;\; \mbox{for}\;0 < \alpha < 1,
\end{eqnarray}
with $\left\{\frac{|x|}{(\eta t^\beta)^{1/\alpha}}\right\} < < 1.$

Further, if  we calculate the residues at the poles of $\Gamma(s/\alpha)$  at the points 
$s = -\alpha\nu\;(\nu =0,1,2,\ldots)$, it gives
\begin{equation}
N(x,t)=\frac{1}{|x|}\sum^\infty_{\nu=0}\frac{\Gamma(1+\alpha\nu)}{\Gamma(1+\beta\nu)\Gamma[1+\frac{\alpha\nu}{2}]\Gamma(-\frac{\alpha\nu}{2})}\left[-\frac{\eta t^\beta}{|x|^\alpha}\right]^\nu,
\end{equation}
for $|x| > 1.$ From (69), it can be readily seen that 
\begin{equation}
N(x,t)\sim \frac{1}{|x|} \;\;\mbox{for large} \;|x|.
\end{equation}

\section{Conclusions}
Reaction-diffusion equations are modeling tools for the dynamics presented by a competition between two or more species, activators and inhibitors or production and destruction, that diffuse in a physical medium. When the diffusion, a homogenizing process, presents characteristic length or time scales that are different for the species, a morphological instability may trigger spatio-temporal pattern formation. This mechanism and its various versions is the Turing instability (Nicolis and Prigogine, 1977; Haken, 2004). The complexity of such equations and systems of equations limits the possibilities to derive analytic closed-form representations for their solutions (Wilhelmsson and Lazzaro, 2001; Kulsrud, 2005; Hundsdorfer and Verwer 2003).

This paper attempts to derive analytic closed-form representations of solutions for a general class of fractional reaction-diffusion equations (eqs. (37), (43), (48)), which also provide as special cases respective solutions of standard reaction-diffusion equations. For this purpose the paper summarizes specific techniques for Laplace, Fourier, and Mellin transforms, results for Mittag-Leffler and Fox functions, as well as the applications of Riemann-Liouville, Weyl, and Caputo fractional calculus for tackling fractional reaction-diffusion equations. Closed-form solutions of the equations are given in terms of Fox's function and their behavior for small and large values of the respective parameter is derived. Fractional order moments and behavior of the fundamental solution as given in in eq. (49) are provided in detail.

The methods used in this paper may also be applied to nonlinear Fokker-Planck equations, taking into account the different physical meaning of reaction, drift, and diffusion coefficients, respectively (Tsallis, 2004; Tsallis and Bukman, 1996).   
\bigskip
\noindent
\section*{References}\par
\medskip
\noindent
Caputo, M.: 1969, \emph {Elasticita  e Dissipazione}, Zanichelli, Bologna.\\ 
Compte, A.: 1996, Stochastic foundations of fractional dynamics, 

\emph {Physical Review E} {\bf 53}, 4191-4193.\\
Del-Castillo-Negrete, D., Carreras, B.A., and Lynch, V.E.: 2003, Front 

dynamics in reaction-diffusion systems with L\'evy flights: A fractional 

diffusion approach, \emph {Physical Review Letters} {\bf 91}, 018302.\\
Del-Castillo-Negrete, D., Carreras, B.A., and Lynch, V.E.: 2002, Front

propogation and segregation in a reation-diffusion model with 

cross-diffusion, \emph {Physica D} {\bf 168-169}, 45-60.\\
Debnath, L.: 2003, Fractional integral and fractional differential equations 

in fluid mechanics, \emph {Fractional Calculus and Applied Analysis} {\bf 6}, 

119-155.\\
Dzherbashyan, M.M: 1966, \emph {Integral Transforms and Representation of} 

\emph {Functions in Complex Domain} (in Russian), Nauka, Moscow.\\
Dzherbashyan, M.M.: 1993, \emph {Harmonic Analysis and Boundary Value} 

\emph {Problems in the Complex Domain}, Birkhaeuser-Verlag, Basel.\\ 
Erd\'elyi, A., Magnus, W., Oberhettinger, F., and Tricomi, F.G.: 1953, \emph {Higher} 

\emph {Transcendental Functions}, Vol. {\bf 1}, McGraw-Hill, New York, Toronto, 

and London.\\
Erd\'elyi, A., Magnus, W., Oberhettinger, F., and Tricomi, F.G: 1953, 

\emph {Higher Transcendental Functions}, Vol. {\bf 2}, McGraw-Hill, New York, Toronto, 

and London.\\
Erd\'elyi, A., Magnus, W., Oberhettinger, F., and Tricomi, F.G: 1954, \emph {Tables} 

\emph {of Integral Transforms}, Vol. {\bf 1}, McGraw-Hill, New York, Toronto, and 

London.\\
Erd\'elyi, A., Magnus, W., Oberhettinger, F., and Tricomi, F.G: 1954b, \emph {Tables} 

\emph {of Integral Transforms}, Vol. {\bf 2}, McGraw-Hill, New York, Toronto, and 

London.\\
Erd\'elyi, A., Magnus, W., Oberhettinger, F., and Tricomi, F.G: 1955, \emph {Higher} 

\emph {Transcendental Functions}, Vol. {\bf 3}, McGraw-Hill, New York, Toronto, 

and London.\\
Haken, H.: 2004, \emph {Synergetics: Introduction and Advanced Topics}, 

Springer-Verlag, Berlin-Heidelberg.\\
Haubold, H.J. and Mathai, A.M.: 1995, A heuristic remark on the periodic 

variation in the number of solar neutrinos detected on Earth,

\emph {Astrophysics and Space Science} {\bf 228}, 113-134.\\ 
Haubold, H.J. and Mathai, A.M.: 2000, The fractional kinetic equation and 

thermonuclear functions, \emph {Astrophysics and Space Science} {\bf 273}, 53-63.\\
Henry, B.I. and Wearne, S.L.: 2000, Fractional reaction-diffusion, \emph {Physica A} 

{\bf 276}, 448-455.\\
Henry, B.I. and Wearne, S.L.: 2002, Existence of Turing instabilities in 

a two-species fractional reaction-diffusion system, 

\emph {SIAM Journal of Applied Mathematics} {\bf 62}, 870-887.\\
Henry, B.I., Langlands, T.A.M., and Wearne, S.L.: 2005, Turing pattern 

formation in fractional activator-inhibitor systems, 

\emph {Physical Review E} {\bf 72}, 026101.\\
Hundsdorfer, W. and Verwer, J.G.: 2003, \emph {Numerical Solution of}

\emph {Time-Dependent Advection-Diffusion-Reaction Equations}, 

Springer-Verlag, Berlin-Heidelberg-New York.\\
Jespersen, S., Metzler, R., and Fogedby, H.C: 1999,  L\'evy flights in external 

force fields: Langevin and fractional Fokker-Planck equations and their 

solutions, \emph {Physical Review  E} {\bf 59}, 2736-2745.\\
Kilbas, A.A and Saigo, M: 2004, \emph {H-Transforms: Theory and Applications}, 

Chapman and Hall/CRC, New York.\\
Kulsrud, R.M.: 2005, \emph {Plasma Physics for Astrophysics}, Princeton 

University Press, Princeton and Oxford.\\
Kuramoto, Y.: 2003, \emph {Chemical Oscillations, Waves, and Turbulence}, 

Dover Publications, Inc., Mineola, New York.\\
Mathai, A.M.: 1993, \emph {A Handbook of Generalized Special Functions for} 

\emph {Statistical and Physical Sciences}, Clarendon Press, Oxford.\\
Mathai, A.M. and Saxena, R.K.: 1978, \emph {The H-function with Applications in} 

\emph {Statistics and Other Disciplines}, Halsted Press [John Wiley and Sons], 

New York, London, and Sydney.\\
Metzler, R. and Klafter, J: 2000, The random walk's guide to anomalous 

diffusion: A fractional dynamics approach, \emph {Physics Reports} 

{\bf 339} 1-77.\\
Metzler, R. and Klafter, J: 2004, The restaurant at the end of the random 

walk: Recent developments in the description of anomalous transport by 

fractional dynamics, \emph {Journal of Physics A: Math. Gen.} {\bf 37}, R161-208.\\
Miller, K.S. and Ross, B.: 1993, \emph {An Introduction to the Fractional Calculus} 

\emph {and Fractional Differential Equations}, John Wiley and Sons, New York.\\
Mittag-Leffler, M.G.: 1903, Sur la nouvelle fonction  $E_\alpha(x)$, 

\emph {Comptes Rendus Acad. Sci. Paris (Ser.II)} {\bf 137}, 554-558.\\
Mittag-Leffler, M.G.: 1905, Sur la representation analytique d'une branche

uniforme d'une fonction monogene, \emph {Acta Mathematica} {\bf 29}, 101-181.\\
Murray, J.D.: 2003, Mathematical Biology, Springer-Verlag, New York.\\
Nicolis, G. and Prigogine, I.: 1977, \emph {Self-Organization in Nonequilibrium}

\emph {Systems: From Dissipative Structures to Order Through Fluctuations}, 

John Wiley and Sons, New York.\\  
Oldham, K.B. and Spanier, J.: 1974, \emph {The Fractional Calculus: Theory and} 

\emph {Applications of Differentiation and Integration to Arbitrary Order}, 

Academic Press, New York; and Dover Publications, New York 2006.\\
Prudnikov, A.P., Brychkov, Yu.A., and Marichev, O.I.: 1989, \emph {Integrals and} 

\emph {Series}, Vol. {\bf 3}, \emph {More Special Functions}, Gordon and Breach, New York.\\
Samko, S.G., Kilbas, A.A., and Marichev, O.I.: 1990, \emph {Fractional Integrals} 

\emph {and Derivatives: Theory and Applications}, Gordon and Breach,

New York.\\ 
Saxena, R.K., Mathai, A.M., and Haubold, H.J.: 2002, On fractional kinetic 

equations, \emph {Astrophysics and Space Science} {\bf 282}, 281-287.\\ 
Saxena, R.K., Mathai, A.M., and Haubold, H.J.: 2004a, On generalized

fractional kinetic equations, \emph {Physica A} {\bf 344}, 657-664.\\
Saxena, R.K., Mathai, A.M., and Haubold, H.J.: 2004b, Unified fractional 

kinetic equation and a fractional diffusion equation, \emph {Astrophysics} 

\emph {and Space Science} {\bf 290}, 299-310.\\
Tsallis, C.: 2004, What should a statistical mechanics satisfy to reflect 

nature?, \emph {Physica D} {\bf 193}, 3-34.\\
Tsallis, C. and Bukman, D.J.: 1996, Anomalous diffusion in the presence 

of external forces: Exact time-dependent solutions and their 

thermostatistical basis, \emph {Physical Review E} {\bf 54}, R2197-R2200.\\
Wilhelmsson, H. and Lazzaro, E.: 2001, \emph {Reaction-Diffusion Problems}

\emph {in the Physics of Hot Plasmas}, Institute of Physics Publishing, 

Bristol and Philadelphia.\\
Wiman, A.: 1905, Ueber den Fundamentalsatz in der Theorie der Functionen

$E_\alpha(x)$, \emph {Acta Mathematica} {\bf 29}, 191-201.\\
Wright, E.M.: 1933, On the coefficients of power series having exponential 

singularities, \emph {Journal of the London Mathematical Society} {\bf 8}, 71-79.\\
Wright, E.M.: 1935, The asymptotic expansion of the generalized hypergeo-

metric functions, \emph {Journal of the London Mathematical Society} 

{\bf 10}, 387-293.\\    
Wright, E.M.: 1940, The asymptotic expansion of the generalized hypergeo-

metric functions, \emph {Proceedings of the London Mathematical Society} 

{\bf 46}, 389-408. 
\end{document}